%% file: 01counterVfin.tex
\documentclass[letterpaper]{amsart}
\usepackage{amsmath,amsthm,amssymb}
\usepackage{tabularx,epsfig,graphicx}
\usepackage[vcentermath]{youngtab}

\begin{document}

\renewcommand{\labelenumi}{\theenumi.}

\bibliographystyle{amsplain}

\input{macros.tex}

\title{Counterexamples to the 0-1 Conjecture}


\begin{abstract}
  For permutations $x$ and $w$, let $\mu(x,w)$ be the coefficient of
  highest possible degree in the Kazhdan-Lusztig polynomial $\pxw$.
  It is well-known that the $\mu(x,w)$ arise as the edge labels of
  certain graphs encoding the representations of $\fsn$.  The 0-1
  Conjecture states that the $\mu(x,w) \in \{0,1\}$.  We present two
  counterexamples to this conjecture, the first in $\fsp{16}$, for
  which $x$ and $w$ are in the same left cell, and the second in
  $\fsp{10}$.  The proof of the counterexample in $\fsp{16}$ relies on
  computer calculations.
\end{abstract}

\author{Timothy J. McLarnan} \email{timm@earlham.edu}

\address{Author's address:
  Dept.\ of Mathematics\\
  Earlham College\\ Richmond, IN 47374}

\author{Gregory S. Warrington} \email{warrington@math.umass.edu}

\address{Author's address:
  Dept.\ of Mathematics \& Statistics\\
  University of Massachusetts\\ Amherst, MA 01003}


\maketitle

\section{Introduction}
\label{sec:intro}

In studying the representations of Hecke algebras, Kazhdan and Lusztig
\cite{K-L1} defined a class of polynomials now known as
Kazhdan-Lusztig (KL) polynomials $\pxw$ that are indexed by pairs of
elements in a Coxeter group.  These polynomials carry important
representation-theoretic and geometric information.  Certain
coefficients $\mu(x,w)$ are particularly important
representation-theoretically in addition to controlling the recursive
structure of the polynomials.  While these coefficients $\mu(x,w)$ are
easily seen to take varied (non-negative) values in most Weyl groups,
empirical evidence has suggested the following
\newtheorem*{zoconj}{0-1 Conjecture}
\begin{zoconj}\label{conj:zerone}
  For $x,w\in\fsn$, $\mu(x,w) \in \{0,1\}$.
\end{zoconj}
If this conjecture were true, Kazhdan and Lusztig's construction of
the irreducible representations of $\fsn$ (see \cite{K-L1}) would be
embodied simply by graphs rather than \textit{edge-labeled} (by the
$\mu(x,w)$) graphs.  However, as the following theorem shows, this is
not the case.
\begin{thm}\label{thm:main}
  Identify elements of $\fsp{16}$ (resp., $\fsp{10}$) with
  permutations of the set $\{0,1,\dotcoms 9,a,b,c,d,e,f\}$ (resp.,
  $\{0,1,\dotcoms 9\}$).  We have the following two equalities:
  \begin{enumerate}
  \item $\mu({\rm 54109832dc76bafe},{\rm c810d942fa53b6e7}) = 5$.
  \item $\mu(4321098765,9467182350) = 4$.
  \end{enumerate}
\end{thm}
The first case offers the smallest counterexample to the 0-1
Conjecture with both permutations lying in the same left cell.  The
existence of such an example implies that the graphs describing the
irreducible representations do, in fact, need to be edge-labeled.
Exhaustive computer calculations by du Cloux \cite{ducloux} and both
authors independently have shown that there are no counterexamples in
$\fsp{9}$ or below.  Hence, the counterexample given in part 2 of
\thmref{thm:main} occurs in the smallest possible group.  The
following corollary is immediate from either part of
\thmref{thm:main}:
\begin{cor}
  The 0-1 Conjecture is false.
\end{cor}
So, in this sense, the combinatorics of the symmetric group is not
simpler than that of other Weyl groups.  

The possibility that $\mu(x,w) \in \{0,1\}$ for any $x,w\in\fsn$ was
noticed by Lascoux and Sch\"utzenberger (see \cite{garmcl}),
presumably noticed by Kazhdan and Lusztig, and certainly noticed
independently by many others.  In fact, Lascoux and Sch\"utzenberger
\cite{LS11} showed the 0-1 Conjecture to be true for Grassmannian
permutations $w$.  However, given the difficulty of examining $\fsn$
for $n \geq 9$ empirically, there has not appeared to be a consensus
as to the truth of the conjecture.

Based on the work of Lascoux and Sch\"utzenberger, Garsia and McLarnan
\cite{garmcl} list three progressively weaker conjectures related to
the 0-1 Conjecture.  The Lascoux-Sch\"utzenberger (L-S) graph
has as vertices all members of a left cell.  All pairs $\{x,w\}$ in
which $w$ covers $x$ in the left weak Bruhat-Chevalley order are edges
of this graph, as are the pairs produced by all possible applications
of the $L_i$ of \defref{def:li}.  Each such edge $\{x,w\}$ is easily
checked to have $\mu(x,w) = 1$.  There are three natural questions
motivated by this construction:
\begin{enumerate}
\item Is this L-S graph identical to the ``K-L graph'' described
  by Kazhdan and Lusztig in \cite{K-L1}?
\item If one starts with the L-S graph and follows the recipe of
  Kazhdan and Lusztig for using the graph to associate a transition
  matrix to each permutation, does one obtain an irreducible
  representation for $\fsn$ corresponding to that left cell?
\item If not, then does one at least get some representation of $\fsn$?
\end{enumerate}
That the first conjecture is false follows from \thmref{thm:main}.1,
since every edge in the L-S graph has weight 1.  In
\secref{sec:norepn}, we sketch computer calculations showing that the
second and third of these conjectures are also false.

\secref{sec:prelim} presents preliminary notation and definitions.
\secref{sec:compproof} describes the algorithm used by the first
author in 1989 to prove \thmref{thm:main}.1.  In
\secref{sec:combproof}, we give the second author's combinatorial
proof of \thmref{thm:main}.2.

\section{Preliminaries}
\label{sec:prelim}

We will consider elements of $\fsn$ as permutations on the set
$\{0,\dotcoms n-1\}$.  As a generating set $\simref$, we will take the
adjacent transpositions $s_i = (i,i+1)$ for $0\leq i \leq n-2$.  A
one-line notation for a permutation $w$ is afforded by writing the
image of the $n$-tuple $[0,\dotcoms n-1]$ under the action of $w$:
$[w(0),w(1),\dotcoms w(n-1)]$ (we often omit the commas and brackets).
The \textit{length function} for $\fsn$ is given by
\begin{equation*}
  l(w) = |\{0 \leq i < j < n: w(i) > w(j)\}|.
\end{equation*}

In \defref{def:bdf}, we define the Bruhat-Chevalley partial order on
$\fsn$.  (The definition we give is equivalent to more common
descriptions such as the tableau criterion; see \cite{BLak,Fulton-book,Hum}
and the references cited therein.)

\begin{dfn}\label{def:rnkdef}
  Let $x,w \in \fsn$, $p,q\in \bbZ$.  Define the \emph{rank function}
  $r_w(p,q) \defeq |\{i \leq p: w(i) \geq q\}|$ and the
  \emph{difference function} $\dxw(p,q) \defeq r_w(p,q) - r_x(p,q)$.
\end{dfn}

\begin{dfn}\label{def:bdf}
  We define the \textit{Bruhat-Chevalley partial order} ``$\leq$'' on
  $\fsn$ by setting $x\leq w$ if and only if $\dxw(p,q) \geq 0$ for
  all $p,q$.
\end{dfn}

Let $[x,w]=\{z: x\le z\le w\}$ in the Bruhat-Chevalley ordering.
Billey and Warrington \cite{gwsb-msl} prove the following result.
\begin{lem}\label{lem:trans}
  If $x(i) = w(i)$, $\dxw(i,x(i)) = 0$ and $z\in [x,w]$, then $z(i)=x(i)$.
\end{lem}

We can view the Bruhat-Chevalley order graphically using ``Bruhat
pictures'' determined by the function $\dxw$.  A typical picture is
shown in \figref{fig:brupic}.  Let $\digw$ refer to the permutation
matrix for $w$.  Entries of $\digx$ (resp., $\digw$) are denoted by
black disks (resp., open circles).  Shading denotes regions in which
$\dxw \geq 1$.  Successively darker shading denotes successively
higher values of $\dxw$.  Positions corresponding to 1's of both
$\digx$ and $\digw$ (termed ``capitols'') are denoted by a black disk
and a larger concentric circle.
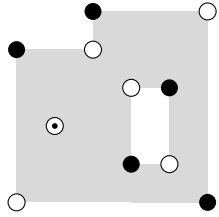
 \mymslfig{brupic}{Bruhat picture for $x = [2,0,4,1,3,5]$, $w =
  [5,2,3,1,4,0]$.}

While not strictly necessary, these pictures help motivate results
such as \lemref{lem:trans} and can be very helpful in computing with
KL polynomials.  In fact, a number of the arguments in
\secref{sec:compproof} were arrived at with the aid of these pictures.
The reader of that section may benefit from constructing the
appropriate Bruhat pictures.

Using the Bruhat-Chevalley order, there are several sets we can
associate to any permutation $w$.  We define the right and left
descent sets of $w$ and the set of flush elements of $w$ to be
\begin{align}
  \rds(w) &= \{s\in\simref: ws < w\},\\
  \lds(w) &= \{s\in\simref: sw < w\} \text{ and }\\
  \ext(w) &= \{x\leq w: \rds(x) \supseteq \rds(w) \text{ and }
                       \lds(x) \supseteq \lds(w)\}.
\end{align}

We now give a combinatorial definition of the Kazhdan-Lusztig (KL)
polynomials applicable to any Coxeter group.  For motivation and a
more natural definition, we refer the reader to \cite{Hum,K-L1}.  In
order to give the definition succinctly, we set
\begin{equation}\label{eq:mudef}
  \mu(x,w) = \text{ coefficient of }
             q^{(l(w)-l(x)-1)/2}\text{ in }\pxw,
\end{equation}
and define $c_s(x) = 1$ if $xs < x$; $c_s(x) = 0$ if $xs > x$. 
\begin{thm}[\cite{K-L1}]\label{thm:kldef}
  There is a unique set of polynomials $\{\pxw\}_{x,w\in\fsn}$ such
    that, for all $x,w\in\fsn$:
  \begin{enumerate}
  \item $\pww = 1$,
  \item $\pxw = 0$ when $x\not\leq w$,
  \item If $s\in\rds(w)$, then
    \begin{equation}\label{eq:std}
      \pxw = q^{c_s(x)}\pxws + q^{1-c_s(x)}\pxsws - \sumsb{z \leq
        ws\\ zs < z} \mu(z,ws)q^{\frac{l(w)-l(z)}{2}} \pxz.      
    \end{equation}
    The analogous recursion with $s$ acting on the left holds when 
    $s\in\lds(w)$.
  \end{enumerate}
  Further, these polynomials satisfy the degree restriction 
  \begin{equation}\label{eq:KLdegree}
    \deg(P_{x,w}) \leq (l(w)-l(x)-1)/2 \text{ when } x < w.
  \end{equation}
\end{thm}
Note that $\mu(x,w)$ is the coefficient of the highest possible power
of $q$ in $\pxw$ and that $\mu(x,w)=0$ if $l(w)-l(x)$ is even.

The complexity of the KL polynomials arises from the sum subtracted
off in \eqref{eq:std}.  We now introduce some notation to let us deal
with these sums concisely.  For $x,w\in\fsn$ and $s\in\rds(w)$, let
\begin{align*}
  \mslor{x}{ws}{s} &= \{z: x \leq z < ws,\ zs < z,\ l(z) < l(ws)-1, 
             z\in\ext(ws)\},\\
  \mscor{x}{ws}{s} &= \{z: x \leq z < ws,\ zs < z,\ l(z) = l(ws)-1\},\\
  \msset{x}{ws}{s} &= \mscor{x}{ws}{s} \cup \mslor{x}{ws}{s},\\
  \musumr{x}{ws}{s} &= \
        \sumsb{z \in \msset{x}{ws}{s}}\mu(z,ws)q^{\frac{l(w)-l(z)}{2}}\pxz.
\end{align*}

\propref{prop:kl}.\ref{item:ext} will imply that $\musumr{x}{ws}{s}$
is the sum appearing in \eqref{eq:std}.  Let $z\in [x,ws]$.  $z$ is
\textit{right $s$-flush} for this interval if $z\in\mslor{x}{ws}{s}$.
It is \textit{right $s$-coatomic} for this interval if
$z\in\mscor{x}{ws}{s}$.  The ``left'' versions are defined analogously
(with ``$(s \cdot)$'' substituted for ``$(\cdot s)$'').  We will omit
``left'' and ``right,'' as they will be clear from context.

We will need several additional properties of KL polynomials that are
not immediately apparent from the definition; we require the following
notation:

\begin{dfn}\label{dfn:fl}
  For $w\in\fsn$ and $0 \leq i_1 < \cdots < i_k \leq n-1$ for $k\leq
  n$, let $\fl[w(i_1),w(i_2),\ldots,w(i_k)]$ be the unique
  \textit{flattened} permutation $[v(1),\ldots, v(k)] \in \fsp{k}$ such
  that $v(j) < v(k)$ precisely when $w(i_j) < w(i_k)$.
\end{dfn}
For example, if $w=7461098253$, then
\begin{equation*}
  \fl[w(0),w(2),w(3),w(5),w(8),w(9)]=\fl[7,6,1,9,5,3]=430521.
\end{equation*}

\begin{dfn}\label{dfn:tilde}
  Let
  \begin{equation*}
    \delxw = \{i: x(i) \neq w(i) \text{ or } \dxw(i,x(i)) \neq 0\}.
  \end{equation*}
  If $\delxw = \{d_1,d_2,\dotcoms d_k\}$ with $d_i < d_j$ for $i <
  j$, we get two \textit{reduced} permutations by flattening $x$ and $w$ with
  respect to $\delxw$:
  \begin{align*}
    \xti &= \fl([x(d_1),x(d_2),\dotcoms x(d_k)]) \ \text{ and }\\
    \wti &= \fl([w(d_1),w(d_2),\dotcoms w(d_k)]).
  \end{align*}
\end{dfn}
Note that $\xti$ and $\wti$ are permutations in $\fsp{k}$.  For
instance, if
\begin{align}\label{eq:unfla}
  x&=6491082753\\
  w&=9461782350,
\end{align}
then $\delxw=\{0,2,4,6,7,9\}$, and
\begin{align}\label{eq:nowfla}
  \xti&=350142\\
  \wti&=534120.
\end{align}
To obtain the Bruhat picture of the pair $\xti$,$\wti$ from that for
$x$,$w$, one simply removes the capitols not surrounded by a
shaded region (see \figref{fig:sampfla}).
\mysmmslfig{sampfla}{Sample Bruhat pictures for \eqref{eq:unfla} and
  \eqref{eq:nowfla}.}

\begin{prop}\label{prop:kl}
  The KL polynomials satisfy the following properties:
  \begin{enumerate}
  \item If $s\in\rds(w)$, then $\pxw = \pxsw$.
        If $s\in\lds(w)$, then $\pxw = \psxw$.\label{item:pxwsim}
  \item $\pxw = \pxiwi$.\label{item:inv}        
  \item $\pxw = \pxtwt$.\label{item:fla}
  \item If $x\not\in\ext(w)$ and $l(x) < l(w)-1$, then 
        $\mu(x,w) = 0$.\label{item:ext}
  \end{enumerate}
\end{prop}
The first two properties are standard and can be found in \cite{K-L1}.
Proof of the third can be found in \cite{gwsb-msl}; the fourth
follows from the first property along with \eqref{eq:KLdegree}.  

The $\mu(x,w)$ also satisfy an identity which we will be integral to
the proof in \secref{sec:compproof}.  To state it, we make the
following definitions:

\begin{dfn}\label{def:li}
  Let $\mathcal{L}_k$ be the set of permutations $w$ for which $s_kw <
  w$ or $s_{k+1}w < w$, but not both.  Define an operator $L_k$ acting
  on $\mathcal{L}_k$ by setting
  \begin{equation*}
    (L_kw)^{-1}(j) =
    \begin{cases}
      w^{-1}(k+2), & \text{ if } j = k,\\
      w^{-1}(k), & \text{ if } j = k+2,\\
      w^{-1}(j), & \text{ otherwise.}
    \end{cases}
  \end{equation*}
\end{dfn}
In other words, $\mathcal{L}_k$ consists of all permutations in which
$k$, $k+1$, $k+2$ do not appear either in increasing or decreasing
order; and $L_kw$ is obtained from $w$ by interchanging $k$, $k+2$.
For instance, $L_2[3,1,4,0,2] = [3,1,2,0,4]$.  The operator $L_k$ is
called an elementary Knuth transformation.  It is intimately connected
to the Robinson-Schensted correspondence discussed below; for details,
see \cite{Fulton-book,Knuth70,Knuth3}.

\begin{dfn}
  For $x$ and $w$ comparable under the Bruhat-Chevalley order, set
  \begin{equation*}
    \mu[x,w] = 
    \begin{cases}
      \mu(x,w), & \text{ if } x \leq w,\\
      \mu(w,x), & \text{ otherwise.}
    \end{cases}
  \end{equation*}
\end{dfn}

\begin{thm}[\cite{K-L1}]\label{thm:samel}
  If $x,w\in \mathcal{L}_k$, then $\mu[x,w] = \mu[L_kx,L_kw]$.
\end{thm}

\section{Computer Proof of \thmref{thm:main}.1}
\label{sec:compproof}

The example of \thmref{thm:main}.1 was found via a computer search for
a counterexample to the 0-1 Conjecture.  As looking at every pair of
permutations even in $\fsp{10}$ is prohibitively expensive, we will
cut down our search space by searching for a counterexample that is
minimal in some sense.  In particular, we will search for a
counterexample $\{x,w\}$ in $\fsn$ with $x$ and $w$ in the same left
cell which minimizes in order the following parameters.
\begin{enumerate}
\item $n$: i.e., $\mu(u,v) \leq 1$ for all $u,v\in\fsp{n-1}$ in the
  same left cell.
\item $l(w) - l(x)$: i.e., $u,v\in\fsn$ in the same left cell with
  $l(v)-l(u) < l(w)-l(x)$ implies $\mu(u,v) \leq 1$.
\item $l(w)$: i.e., $v\in\fsn$ with $l(v) < l(w)$ implies that there
  does not exist a $u\in\fsn$ in the same left cell as $v$ having
  $\mu(u,v) > 1$.
\end{enumerate}

The key to searching efficiently turns out to be the
Robinson-Schensted correspondence, which we now recall.  The material
in this section will be presented briefly---a more detailed
exposition can be found in \cite{garmcl} (also see \cite{Fulton-book,Knuth3}).  

Let $\lambda$ be a partition of $n$ (denoted $\lambda \vdash n$) with
$\lambda_1 \geq \lambda_2 \geq \cdots \lambda_k \geq 0$.  We associate
a Ferrers diagram consisting of left-justified rows of boxes with
$\lambda_i$ boxes in the $i$-th row from the bottom.  A
\textit{standard tableau of shape $\lambda$} is an injective filling
of these boxes with $0,1,2,\dotcoms n-1$ such that the entries
increase from left to right on rows and from bottom to top on columns.
The \textit{column word}, $\cwd(T)$, of a tableau $T$ is obtained by
reading the columns of $T$ from top to bottom starting with the
leftmost column.  The \textit{row word}, $\rwd(T)$, of $T$ is obtained
by reading the rows of $T$ from left to right starting with the top
row.  Any word that can be obtained in this way from some tableau is
called a {\it tableau word}.  The \textit{descent set $D(T)$} of a
tableau $T$ is the set of indices $i$ for which $i+1$ is strictly to
the north and weakly to the west of $i$.  In other words,
\begin{equation*}
  D(T)=D_L(\rwd(T))=D_L(\cwd(T)).
\end{equation*}

The Robinson-Schensted correspondence gives a bijection between the
elements $w\in\fsn$ and the pairs of tableaux of the same shape
$\lambda \vdash n$.  For the specifics of the bijection, see, e.g.,
\cite{Fulton-book}.  Via the Bruhat-Chevalley order on permutations,
this correspondence induces a partial order on pairs of tableaux of
the same shape which we will also denote by ``$\leq$''. The
\textit{left cell} indexed by the tableau $\rtq$ consists of all pairs
$(\ltw,\rtq)$ where $\ltw$ has the same shape as $\rtq$.  Below we
illustrate these definitions:

\begin{align}\label{eq:tabex}
  x = 4265013\leftrightarrow 
  (\ltx,\rtx) &= \left(
    \young(46,25,013)\ ,\ \young(45,13,026)\right),\\
  \rwd(\ltx) &= 4625013,\\
  \cwd(\ltx) &= 4206513.
\end{align}

In order to describe the edge-labeled graphs defined by Kazhdan and
Lusztig corresponding to \textit{irreducible} representations of
$\fsn$, we need only consider pairs $x,w$ lying in the same left cell.
Thus, a counterexample among such pairs shows that the edge-labeling
is necessary.

Fortunately, there is complete redundancy amongst the left cells
with respect to the values of the $\mu[x,w]$:
\begin{thm}[\cite{K-L1}]\label{thm:col}
  Let $\ltx$, $\ltw$, $\rtq$ and $\rtqp$ be tableaux of the same shape.
  Let $x$, $w$, $x'$ and $w'$ correspond, under the Robinson-Schensted
  correspondence, to the following pairs of tableaux:
  \begin{align*}
    x \leftrightarrow (\ltx,\rtq),\, & \qquad \ w \leftrightarrow (\ltw,\rtq),\\
    x' \leftrightarrow (\ltx,\rtqp), & \qquad w' \leftrightarrow (\ltw,\rtqp).
  \end{align*}
  Then $\mu[x,w] = \mu[x',w']$.
\end{thm}

With \thmref{thm:col}, to search for a minimal counterexample, we
effectively need only search over pairs $(\ltx,\rtq),(\ltw,\rtq)$ of
tableaux of the same shape without regard to $\rtq$.  One can get a
sense of the savings by noting that $|\fsp{16}| = 20,922,789,888,000$,
while the number of standard Young tableau of size $16$ is a mere
$46,206,736$.  The task of considering all \textit{pairs} of
permutations is clearly infeasible.  While one still needs to look at
pairs of tableaux (of the same shape), even naively we need only
consider $|\fsp{16}|$ pairs.  And, with the proper filters, we can do
much better.  

We are now ready to present the main facts upon which the algorithm
rests.
\begin{lem}\label{lem:eprops}
  Let $x < w \in \fsn$ with 
  \begin{equation*}
    x \leftrightarrow (\ltx,\rtq)\,  \quad \text{ and }\quad \ w \leftrightarrow (\ltw,\rtq)
  \end{equation*}
  under the Robinson-Schensted correspondence.  If the pair $\{x,w\}$ is a
  counterexample to the 0-1 Conjecture satisfying the above three
  minimality properties, then the following eight conditions must hold:
  \begin{enumerate}
  \item $D(\ltw) \subseteq D(\ltx)$.\label{item:p1}
  \item $(\ltw, \rtqp) > (\ltx, \rtqp)$ for all tableaux $\rtqp$.\label{item:p2}
  \item The largest number, $n - 1$, sits strictly higher in $\ltw$
    than in $\ltx$.\label{item:p3}
  \item If $\rwd(\ltw)^{-1}(k+2) < \rwd(\ltw)^{-1}(k)$, then\label{item:p4}
    \begin{equation*}
      \rwd(\ltx)^{-1}(k+2) < \rwd(\ltx)^{-1}(k+1) < \rwd(\ltx)^{-1}(k).
    \end{equation*}
  \item There do not exist $L_{i_1}, L_{i_2} , \dotcoms L_{i_k}$ such
    that that\label{item:p5}
    \begin{equation}\label{eq:lineq}
      l(L_{i_k}\cdot \cdots L_{i_2} L_{i_1}\ltw) - 
      l(L_{i_k}\cdot \cdots L_{i_2} L_{i_1}\ltx) < l(\ltw) - l(\ltx). 
    \end{equation}  
  \item There do not exist $L_{i_1}, L_{i_2} , \dotcoms L_{i_k}$
    satisfying both \eqref{eq:lineq} and \label{item:p6}
  \begin{equation}
    l(L_{i_j} L_{i_{j-1}} \cdots L_{i_1} \ltw) - l(L_{i_j} L_{i_{j-1}}
    \cdots L_{i_1} \ltx) = l(\ltw) - l(\ltx), \text{ all } j < k.
  \end{equation}
  \item $l(w) - l(x)$ is odd.     \label{item:p7}
  \item For no $i$ are $0, 1,\dotcoms i - 1$ in identical positions in
    $w$ and $x$ and are $\wti$ and $\xti$ tableau words of the same
    shape.  \label{item:p8}
  \end{enumerate}
\end{lem}

\begin{proof}
  We give a brief justification for each condition:
  \begin{enumerate}
  \item This follows from \propref{prop:kl}.\ref{item:ext} and the
    fact that $D(\ltw) = D_L(\rwd(\ltw))$.

  \item Incomparability for some $\rtqp$ would imply a contradiction
    by \thmref{thm:col}.  So consider the case where $(\ltw,\rtqp) <
    (\ltx,\rtqp)$ for some $\rtqp$.  As detailed in \cite[Section
    5]{garmcl}, this implies that there exists some $\rtqpp$ for which
    $w'' \leftrightarrow (\ltw,\rtqpp)$ and $x'' \leftrightarrow
    (\ltw,\rtqpp)$ are related in the \textit{weak} Bruhat-Chevalley order and
    satisfy $l(w'') - l(x'') = 1$.  But then, by definition,
    $\mu(x'',w'') = 1$.  A contradiction then results by applying
    \thmref{thm:col}.

  \item By the previous property, we must have $\cwd (\ltw) \geq \cwd
    (\ltx)$.  This implies that $n-1$ must be at least as high in
    $\ltw$.  If it is the same height, then by
    \propref{prop:kl}.\ref{item:fla}, you could delete it and get a
    counterexample in $\fsp{n-1}$.

  \item Given minimality property 3, this is equivalent to the
    first property along with \thmref{thm:col}.

  \item Knuth transformations preserve left cells; so by
    \thmref{thm:col}, existence would contradict minimality property
    2.

  \item This is a special case of the previous property.

  \item If $l(w) - l(x)$ is even, then $\mu(x,w) = 0$.

  \item Otherwise, $\xti,\wti$ lie in the same left cell and afford a
    smaller counterexample by \propref{prop:kl}.\ref{item:fla}. (Note
    that $l(w)-l(x) = l(\wti)-l(\xti)$.)
  \end{enumerate}
\end{proof}

We were able to write code to check quickly whether a pair $(x,w)$
satisfies Properties \ref{item:p1}, \ref{item:p3}, \ref{item:p4},
\ref{item:p7} and \ref{item:p8}. It's more time-consuming to check
\ref{item:p6}.  We found it slowest to check Properties \ref{item:p5}
and \ref{item:p2}.

Since we are working with such large groups, considerable care must be
taken to check each of the above eight properties as efficiently as
possible.  For instance, it is impossibly slow to test property
\ref{item:p2} by computing all $\rtqp$ tableaux, doing inverse
Robinson-Schensted, and checking the Bruhat-Chevalley relations.  It
is much faster to generate the pairs by doing Knuth transformations
and to check for the Bruhat-Chevalley relation by seeing whether the
Knuth transformation has destroyed the Bruhat-Chevalley relation which
applied before the transformation.

The algorithm used to find a counterexample is as follows:
\begin{enumerate}\renewcommand{\labelenumi}{Step \theenumi.}
\item Build up pairs of tableau $\ltx$ and $\ltw$ that satisfy
  properties \ref{item:p1}, \ref{item:p3}, \ref{item:p4} and
  \ref{item:p8} one letter at a time.
\item Successively filter out those pairs not satisfying each of
  properties \ref{item:p7}, \ref{item:p6}, \ref{item:p5} and
  \ref{item:p2}.
\item For all remaining pairs $\ltx$ and $\ltw$, choose $\rtqp$ to
  minimize the length difference between $x \leftrightarrow (\ltx,\rtqp)$
  and $w \leftrightarrow (\ltw,\rtqp)$.
\item Compute $\mu(x,w)$.  Filter out those pairs for which $\mu(x,w)
  \leq 1$.
\end{enumerate}
No pairs in $\fsp{13}$ or below make it through Step 2.  In
$\fsp{14}$ and $\fsp{15}$, none make it through Step 4.  But in
$\fsp{16}$, the following pair of permutations passes all steps:
\begin{alignat*}{2}
    w &= {\rm c810d942fa53b6e7}, & \quad l(w) = 53,\\
    x &= {\rm 54109832dc76bafe}, & \quad l(x) = 32.
\end{alignat*}
The difference in lengths is 21, and the leading coefficient (the
coefficient of degree 10) of $\pxw(q)$ is $\mu(x, w) = 5$. The K-L
polynomial in its entirety is
\begin{equation*}
\pxw(q) = 5q^{10} +72q^9 + 387q^8 + 1039q^7 +1610q^6 +1536q^5 + 931q^4 + 365q^3 +
92q^2 +14q + 1. 
\end{equation*}
This completes the proof of \thmref{thm:main}.1.

Exactly one other pair of permutations in $\fsp{16}$ passes all the
steps of this algorithm, affording a second counterexample to the 0-1
Conjecture:
\begin{alignat*}{2}
    w &= {\rm ca610fb732d84e95}, & \quad l(w) = 60,\\
    x &= {\rm 76310cb542a98fed}, & \quad l(x) = 39,
\end{alignat*}
\begin{equation*}
\pxw(q) = 5q^{10} +56q^9 + 231q^8 + 533q^7 +776q^6 +755q^5 + 501q^4 + 226q^3 +
67q^2 +12q + 1. 
\end{equation*}

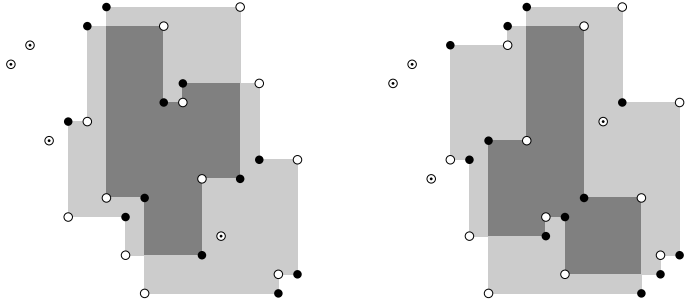
\mytinymslfig{mcl}{The Bruhat pictures for the two minimal
  counterexamples to the 0-1 Conjecture lying in $\fsp{16}$.}

Although we have not been able completely to verify the results in
this section without the use of computers, we are extremely confident
of the truth of \thmref{thm:main}.1.  The two authors began
collaborating after we had independently written programs to compute
K-L polynomials, and these programs agree on the values of the
polynomials computed above.  Only the first author has carried out the
process of generating and filtering pairs to produce these
counterexamples, but errors in that code would only affect the
minimality of our examples.  It seems extraordinarily unlikely that
our completely independent computations of the K-L polynomials could
be incorrect and yet agree.  The computer code used in the proof of
\thmref{thm:main}.1 (along with java code for computing K-L
polynomials) is archived in the source package for this paper on
http://arXiv.org.

\section{Combinatorial Proof of \thmref{thm:main}.2}
\label{sec:combproof}

The previous section describes a counterexample showing that for large
enough~$n$, labeled graphs are required for Kazhdan and Lusztig's
description of the irreducible representations of $\fsn$.
In this section, we remove the condition that the two
permutations lie in the same left cell.  Such a counterexample
is less interesting representation-theoretically, but it can be
carried out entirely by hand, and it should lend
insight into how and when $\mu(x,w)$ can be greater than 1.  The
counterexample we present was arrived at by close examination (using
the ``Bruhat pictures'' of \cite{gwsb-msl}) of the counterexample
presented in the previous section.  

We begin our proof of \thmref{thm:main}.2 by first calculating several
intermediate KL polynomials.  The main tools in the proof are the
defining recurrence relation \eqref{eq:std} and parts
\ref{item:pxwsim} and \ref{item:fla} of \propref{prop:kl}.  For each
application of \eqref{eq:std}, there are usually several choices of
the generator $s$.  While our choices for $s$ may seem \textit{ad
hoc}, they are actually carefully made both to maximize the number of
applications of \propref{prop:kl}.\ref{item:fla} we can make and to
simplify the calculation of the resulting $\Theta$'s.

\begin{lem}\label{lem:simp}
  The following equalities hold:
  \begin{enumerate}
  \item $\pstack{1032,3120} = 1 + q$.
  \item $\pstack{0213,2301} = 1 + q$.
  \item $\pstack{315042,534120} = 1 + 3q + q^2$.
  \item $\pstack{3106542, 6345120} = 1 + 4q + 4q^2 + q^3$.
  \end{enumerate}
\end{lem}

In the rest of this paper, for layout reasons, we sometimes write
$\pstack{x\\w}$ for $\pxw$.
\begin{proof}
  The first two equalities can be shown immediately using
  \eqref{eq:std} or they can be found in \cite{gwsb-msl}.  For the
  third equality, we begin by expanding using \eqref{eq:std} with $s =
  s_4$:
  \begin{equation}
  \pstack{315042\\534120} = q\pstack{315042\\ 534102} + 
                  \pstack{315024\\ 534102} - 
                  \musumr{315042}{534102}{s_4}.
  \end{equation}  
  Consider $\musumr{315042}{534102}{s_4}$.  It is clear that there are
  no $s_4$-coatomic elements.  By \lemref{lem:trans}, for any 
  $z\in[315042,534102]$, we need $z(5) = 2$.  And for $z$ to be
  $s_4$-flush, we need $\rds(z) \supseteq \{s_0,s_2,s_3,s_4\}$.
  But these conditions cannot simultaneously be satisfied along with
  $z(0) \geq 3$ (necessary for $z \geq 315042$).  So
  $\msset{315024}{534102}{s_4} = \emptyset$ and
  $\musumr{315024}{534102}{s_4} = 0$.  Therefore,
  \begin{align*}
  \pstack{315042\\534120} 
    &= q\pstack{21403\\ 42310} + \pstack{315204\\ 534102} 
        \qquad \qquad \qquad \qquad \quad \ \ \ \ \ \ \ \,
        \text{ (\propref{prop:kl}, parts \ref{item:pxwsim} \&
                                            \ref{item:fla})}\\
    &= q\pstack{21430\\ 42310} + \pstack{20413\\42301} 
        \qquad \qquad \qquad \qquad \qquad \ \ \ \ \ \ 
        \text{ (\propref{prop:kl}, parts \ref{item:pxwsim} \&
                        \ref{item:fla})}\\
    &= q\pstack{1032\\3120} + \left(q\pstack{20413\\24301} +
                                        \pstack{02413\\24301} -
                                        \musumr{20413}{24301}{s_0}\right).\\
  \end{align*}  
  (The last equality follows from \propref{prop:kl}.\ref{item:fla}
  and \eqref{eq:std} with $s = s_0$.)
  
  It is clear that there are no $s_0$-coatomic elements.  If
  $z\in[20413,24301]$, then it follows from \lemref{lem:trans} that
  $z(0) = 2$.  But for $z$ to be $s_0$-flush we need $\rds(z)
  \supseteq \{s_0,s_1,s_2\}$.  These two conditions cannot
  simultaneously be satisfied. Therefore, $\msset{20413}{24301}{s_0} =
  \emptyset$ and $\musumr{20413}{24301}{s_0} = 0$.  Therefore,
  \begin{align*}
  \pstack{315042\\534120} 
    &= q(1+q) + (q\pstack{0312\\3201} + 
                                        \pstack{04213\\ 24301})
        \quad \, \text{ (\lemref{lem:simp}.1;
          \propref{prop:kl}, parts \ref{item:pxwsim} \& \ref{item:fla})}\\
    &= q(1+q) + (q\cdot 1 + \pstack{0213\\ 2301})
        \quad \quad \ \text{ (\thmref{thm:kldef}.1;
          \propref{prop:kl}, parts \ref{item:pxwsim} \&
          \ref{item:fla})}\\
    &= q(1+q) + (q\cdot 1 + (1+q))
        \qquad \qquad \qquad \qquad \qquad \qquad \quad \ \ \ 
        \text{ (\lemref{lem:simp}.2)}\\
    &= 1 + 3q + q^2.     
  \end{align*}
  
  For the fourth equality, we set $s = s_3$ in \eqref{eq:std}, and
  utilize \propref{prop:kl}.\ref{item:pxwsim}:
  \begin{equation*}
    \pstack{3106542\\ 6345120} = (1+q)\pstack{3106542\\ 6341520}
                                - \musumr{3106542}{6341520}{s_3}.
  \end{equation*}
  There are no $s_3$-coatomic elements in $[3106542,6341520]$.  If $z$
  is in this interval, then by \lemref{lem:trans}, $z(4) = 5$.  For
  $z$ to be $s_3$-flush, it must satisfy $\rds(z) \supseteq
  \{s_0,s_2,s_3,s_4,s_5\}$.  As these conditions cannot be
  simultaneously satisfied, we conclude that
  $\msset{3106542}{6341520}{s_3} = \emptyset$; hence
  $\musumr{3106542}{6341520}{s_3} = 0$.  Therefore,
  \begin{align*}
    \pstack{3106542\\ 6345120} 
    &= (1+q)\pstack{3106542\\6341520}\\
    &= (1+q)\pstack{315042\\534120} \
        \qquad \qquad \qquad \qquad \qquad \quad \ \ \ \ 
        \text{ (\propref{prop:kl}, parts \ref{item:pxwsim} \& \ref{item:fla})}\\
    &= (1+q)(1 + 3q + q^2)
        \qquad \qquad \qquad \qquad \qquad \qquad \qquad \quad \ \ \         
        \ \text{ (\lemref{lem:simp}.3)}\\
    &= 1 + 4q + 4q^2 + q^3.     
  \end{align*}
\end{proof}

In addition to the above KL polynomials, we also need to compute
several $\Theta$'s.

\begin{lem}\label{lem:musums}
The following equalities hold:
\begin{enumerate}
\item $\musuml{32170654}{72561340}{s_2} = q^4$.
\item $\musumr{321087654}{835617240}{s_4} = 0$.
\item $\musumr{4321098765}{9461782350}{s_3} = q^4(1+q)$.
\end{enumerate}
\end{lem}

\begin{proof}
  It is easily checked that there are no coatomic elements for any of
  the three above cases.  Hence, in the following, we will assume that
  $z$ is $s$-flush.

\begin{enumerate}
\item Let $x = 32170654$ and $v = s_2w = 72561340$.  To find the
  elements of $\mslol{x}{v}{s_2}$ is straightforward but tedious.  If
  $z\in\mslol{x}{v}{s_2}$, then helpful facts about $z$ include
  \begin{enumerate}
    \item $\{s_0,s_1,s_2,s_4,s_6\}\subseteq \lds(z)$.\label{item:lds}
    \item $\{s_0,s_3,s_6\} \subseteq \rds(z)$.
    \item $z(1)=2$, which combined with \ref{item:lds} implies that
    $z(0)=3$.
    \item $z(2)\in\{1,5,7\}$.
    \item To have $z\in[x,v]$ requires that $z^{-1}(1)\in\{2,3,4\}$,
    that $z^{-1}(7)\in\{2,3\}$, and that $z^{-1}(6)\in\{3,4,5\}$.
  \end{enumerate}
  Armed with these facts, it is not hard to find the nine elements of
  $\mslol{x}{v}{s_2}$.  Of these nine, only three have an odd length
  difference with respect to $v$ (an even length difference with
  respect to $w$); only these three, which are shown in Table
  \ref{tab:threez}, can contribute to
  $\musuml{32170654}{72561340}{s_2}$.

\begin{table}
  \begin{center}
    \begin{tabular}[c]{|>{$}c<{$}>{$}c<{$}|
        >{$}c<{$}>{$}c<{$}|}\hline
& \boldsymbol{z} & \boldsymbol{\zti} & \boldsymbol{\vti}\\\hline
z_1 & 32170654 & 2160543 & 6451230\\
z_2 & 32175640 & 10423 & 42301\\
z_3 & 32751640 & 0312 & 3120\\\hline
    \end{tabular}
  \end{center}
  \caption{Cases for \lemref{lem:musums}.1.}  \label{tab:threez}
\end{table}

We know from \lemref{lem:simp}.4 and \propref{prop:kl}, parts
\ref{item:inv} and \ref{item:fla} that $\pstack{z_1,v} = 1+4q
+4q^2+q^3$.  As $l(v)-l(z_1) = 7$, $\mu(z_1,v) = 1$.  Finally, since
$P_{z_1,z_1} = 1$, the only non-zero term of $\musuml{x}{v}{s_2}$ is
$1 \cdot q^4 \cdot 1 = q^4$.

\item If $z\in[321087654,835617240]$ then by \lemref{lem:trans}, $z(5)
  = 7$.  And if $z\in\mslor{321087654}{835617240}{s_4}$, then
  $\rds(z) \supseteq \{s_0,s_3,s_4,s_5,s_7\}$.  These two conditions
  cannot be satisfied simultaneously.
  
\item Let $x = 4321098765$ and $v = ws_3 = 9461782350$.  To find the elements
  of $\mslor{x}{v}{s_3}$ is again straightforward but still more
  tedious.  If $z\in\mslor{x}{v}{s_3}$, then
\begin{enumerate}
\item $\rds(z) \supseteq \{s_0,s_2,s_3,s_5,s_8\}$.
\item $\lds(z) \supseteq \{s_0,s_3,s_5,s_8\}$.
\item $z(3)=1$; $z(4)=0$.
\item $z^{-1}(8) \leq 6$; $z^{-1}(7)\leq 7$.
\item $z^{-1}(9)\in\{0,2,5\}$.
\item $z^{-1}(8)\in\{5,6\}$.
\item $z(9)\in\{2,3,5\}$.
\end{enumerate}

These observations let us generate the 34 elements of
$\mslor{x}{v}{s_3}$.  Since $l(v) = 30$, we can only have $\mu(z,v)
\neq 0$ if $l(z)$ is odd.  In Table \ref{tab:zlist}, we list the
seventeen of these $z$ with an odd length difference with respect to
$v$ along with the corresponding $\zti$ and $\vti$.

\begin{table}
  \begin{center}
    \begin{tabular}[c]{|>{$}c<{$}>{$}c<{$}|
        >{$}c<{$}>{$}c<{$}|}\hline
& \boldsymbol{z} & \boldsymbol{\zti} & \boldsymbol{\vti}\\\hline
z_1 &    4371098265 & 326087154 & 8356712410\\
z_2 &    4371098652 & 32507641 & 73456120   \\
z_3 &    4391087265 & 32706154 &  73561240  \\
z_4 &    4391087652 & 3260541  & 6345120    \\
z_5 &    6421098753 & 3106542  & 6345120    \\
z_6 &    6471098352 & 230541   & 523410     \\
z_7 &    6491082753 & 350142   & 534120     \\
z_8 &    6491087352 & 24031    & 42310      \\
z_9 &    7431098265 & 52076143 & 74561230   \\
z_{10} & 7431098652 & 4206531  & 6345120    \\
z_{11} & 7461098253 & 305412   & 534120     \\
z_{12} & 9421083765 & 102543   & 451230     \\
z_{13} & 9421086753 & 10342    & 34120      \\
z_{14} & 9431087265 & 205143   & 451230     \\
z_{15} & 9431087652 & 20431    & 34120      \\
z_{16} & 9461083752 & 0231     & 3120       \\
z_{17} & 9461087253 & 0312     & 3120       \\\hline      
    \end{tabular}
  \end{center}
  \caption{Cases for \lemref{lem:musums}.3.}  
\end{table}\label{tab:zlist}

By \propref{prop:kl}.\ref{item:ext}, we ascertain that the only $z$ in
the above table for which we might have $\mu(z,v) \neq 0$ is $z_5 =
6421098753$.  By \propref{prop:kl}.\ref{item:fla} and
\lemref{lem:simp}.4, $\pkl{z_5}{v} = 1 + 4q+4q^2+q^3$.  As
$l(v)-l(z_5) = 7$, $\mu(z_5,v) = 1$.
\propref{prop:kl}.\ref{item:fla}, along with \lemref{lem:simp}.1,
shows that $\pkl{x}{z_5} = 1+q$.  The only non-zero contribution to
the sum in \eqref{eq:std} is therefore $1 \cdot q^4 \cdot (1+q) =
q^4(1+q)$, as desired.
\end{enumerate}

\end{proof}

\begin{proof}[Proof of \thmref{thm:main}.2]
  By \eqref{eq:std} with $s=s_3$ and
  \propref{prop:kl}.\ref{item:pxwsim},
\begin{equation}
  \label{eq:big1}
  \pstack{4321098765\\9467182350} =
(1+q)\pstack{4321098765\\9461782350} - 
\musumr{4321098765}{9461782350}{s_3}.
\end{equation}

By \lemref{lem:musums}.3, $\musumr{4321098765}{9461782350}{s_3} =
q^4(1+q)$.  Using \propref{prop:kl}.\ref{item:fla}, we can therefore
rewrite \eqref{eq:big1} as:
\begin{equation}
    \pstack{4321098765\\9467182350} =
    (1+q)\pstack{321087654\\835671240} - q^4(1+q).
\end{equation}
Expanding using \eqref{eq:std} with $s = s_4$ and applying
\propref{prop:kl}.\ref{item:pxwsim} and \lemref{lem:musums}.2, we get
\begin{align}
  \pstack{4321098765\\9467182350} &=
  (1+q)\left((1+q)\pstack{321087654\\835617240} -
    0\right) - q^4(1+q).\\
  \intertext{By \propref{prop:kl}, parts \ref{item:pxwsim} and
    \ref{item:fla}, this can be rewritten:}
  \pstack{4321098765\\9467182350}
  &= (1+q)(1+q)\pstack{32170654\\73561240} - q^4(1+q)\label{eq:bil1}\\
  &= (1+q)^2\left((1+q)\pstack{32170654\\72561340}-
    q^4\right) - q^4(1+q).\\
  \intertext{The second follows from the first by the left-hand
    version of \eqref{eq:std} with $s=s_2$ and \lemref{lem:musums}.1.
    Simplifying according to \propref{prop:kl}.\ref{item:fla}, we get}
  \pstack{4321098765\\9467182350}
  &= (1+q)^2\left((1+q)\pstack{2160543\\6451230} - q^4\right) - q^4(1+q)\label{eq:bil3}\\
  &= (1+q)^2\left((1+q)(1+4q+4q^2+q^3) - q^4\right) - q^4(1+q)\label{eq:bil4}\\
  &= 1 + 7q + 19q^2 + 26q^3 + 17q^4 + 4q^5.\label{eq:bil5}
\end{align}

Going from \eqref{eq:bil3} to \eqref{eq:bil4} uses
\lemref{lem:simp}.4 and \propref{prop:kl}.\ref{item:inv}.  As
$l(9467182350) - l(4321098765) = 11$, this completes the proof of the
theorem.
\end{proof}

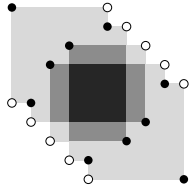
\mytinymslfig{mydual}{The Bruhat picture corresponding to the
  counterexample to the 0-1 Conjecture lying in $\fsp{10}$.}

\begin{remark}
  $\mu(x,w)$ can, in fact, be 2 or 3, though we have not yet found any
  examples in groups smaller than $\fsp{14}$.
\begin{alignat*}{2}
  w &= {\rm 789ab0cd123456}, & \quad l(w) = 47,\\
  x &= {\rm 0759321cba486d}, & \quad l(x) = 32,
\end{alignat*}
\begin{equation*}
\pxw(q) = 2q^7 + 111q^6 + 693q^5 + 1292q^4 + 908q^3 + 257q^2 + 29q + 1.
\end{equation*}

\begin{alignat*}{2}
  w &= {\rm 789ab0cd123456}, & \quad l(w) = 47,\\
  x &= {\rm 0784321cba956d}, & \quad l(x) = 32,
\end{alignat*}
\begin{equation*}
\pxw(q) = 3q^7 + 124q^6 + 716q^5 + 1346q^4 + 960q^3 + 263q^2 + 29q + 1.
\end{equation*}
\end{remark}

\section{Computations on Remaining Conjectures}
\label{sec:norepn}

In \secref{sec:intro}, we mention the conjectures that the L-S graph
obtained by taking the trivial edges in a left cell and adding all the
edges obtained from these by Knuth relations might be the same as the
K-L graph, or at least that this graph might give rise to a
representation of $\fsp{n}$.  That the first of these conjectures is
false follows at once from the counterexample in \thmref{thm:main}.1,
since every edge of the L-S graph has weight 1.  Thus, this conjecture
must be false starting in $\fsp{16}$.  In fact, a computer search for
counterexamples to this conjecture, carried out in a manner analogous
to that in \secref{sec:compproof}, shows that the first counterexamples
to this conjecture appear in $\fsp{14}$.  We obtained these
counterexamples before the counterexamples to the 0-1 Conjecture, and
their existence inspired us to continue searching for edges of weight
greater than 1.

A typical edge in $\fsp{14}$ that cannot be obtained from any trivial
edge by a sequence of Knuth relations is that joining
$w={\rm db630c7418295a}$ and $x={\rm 76530db4192c8a}$.

Not only is the L-S graph not the same as the K-L graph, but the L-S
graph fails to give rise to a representation of $\fsp{n}$, again
starting at $n=14$.  This has
been shown by the first author via computer calculations.  In
particular, take the tableaux whose column words are the permutations
in the counterexample above,
\begin{equation*}
  \ltw = \young(d,bc,67,3489,0125a) \text{ and }
  \qquad \ltx = \young(7,6d,5b,349c,0128a).
\end{equation*}
Choose any tableau $Q$ of the same shape as $\ltx$ and $\ltw$ and let
\begin{equation*}
  x \leftrightarrow (\ltx,\rtq) \text{ and } w \leftrightarrow (\ltw,\rtq)
\end{equation*}
under the Robinson-Schensted correspondence.  For an adjacent
transposition $s$, let $A(s)$ be the matrix obtained by following the
recipe of Kazhdan and Lusztig starting with the L-S graph, and let
$B(s)$ be the matrix obtained by following the
recipe of Kazhdan and Lusztig starting with the K-L graph.  Since the
$B(s)$ generate a representation, we know that for every commuting $s$
and $t$ we have the equality of the matrix entries
\begin{equation}\label{eq:bprod}
  (B(s)B(t))_{w,tx} = (B(t)B(s))_{w,tx}.
\end{equation}

If the $A(s)$ also generate a representation, then we should also have 
\begin{equation}\label{eq:aprod}
  (A(s)A(t))_{w,tx} = (A(t)A(s))_{w,tx}
\end{equation}
for all choices of commuting $s$ and $t$.

If $s = s_4$ and $t = s_b$, then the edge between $x$ and $w$, present
in the K-L graph and absent in the L-S graph, makes a contribution to
the left hand side of \eqref{eq:bprod} via the term
\begin{equation*}
B(s)_{w,x}B(t)_{x,tx}=1\cdot1=1.
\end{equation*}
The corresponding term on the left hand side of \eqref{eq:aprod} is
\begin{equation*}
A(s)_{w,x}A(t)_{x,tx}=0\cdot1=0.
\end{equation*}
It is not obvious that this edge forces different contributions to the
right hand sides of \eqref{eq:bprod} and \eqref{eq:aprod}.  This
suggests that for these values of $x$ and $w$ and $s$ and $t$,
\eqref{eq:aprod} might well be false.

The representation here has dimension 68,640, and finding each matrix
entry involves computing a K-L polynomial; but the matrices are
sufficiently sparse that nearly every term in the matrix products in
\eqref{eq:aprod} is obviously 0.  It is therefore not difficult to
show by computer that for this $x, w, s, t$ we have
\begin{equation*}
  (A(s)A(t))_{w,tx} = 0 \ne 1 = (A(t)A(s))_{w,tx}.
\end{equation*}
Thus, this conjecture, like all those mentioned here, is also false.
  
\section{Acknowledgements}
The first author would like to thank his wife Ann, who shared a
bedroom with many of the computer processor cycles used in this work,
and to thank Adriano Garsia, who started him thinking about K-L
representations.  The second author would like to thank his wife Jill
and also Sara Billey with whom several of the techniques used in this
paper were developed.

\bibliography{gen}

\end{document}

%% file: macros.tex


\hfuzz4pt 

\setlength{\floatsep}{15pt plus 12pt}
\setlength{\textfloatsep}{\floatsep}

\setlength{\intextsep}{\floatsep}\def\ldotsplus{\mathinner{\ldotp\ldotp\ldotp\ldotp}}
\setlength{\intextsep}{\floatsep}\def\ldotscomm{\mathinner{\ldotp\ldotp\ldotp,}}
\def\fourdots{\relax\ifmmode\ldotsplus\else$\m@th \ldotsplus\,$\fi}
\def\dotcoms{\relax\ifmmode\ldotscomm\else$\m@th \ldotsplus\,$\fi}

%
%
\setlength{\unitlength}{0.06em}
\newlength{\cellsize} \setlength{\cellsize}{18\unitlength}
\newsavebox{\cell}
\sbox{\cell}{\begin{picture}(18,18)
\put(0,0){\line(1,0){18}}
\put(0,0){\line(0,1){18}}
\put(18,0){\line(0,1){18}}
\put(0,18){\line(1,0){18}}
\end{picture}}
\newcommand\cellify[1]{\def\thearg{#1}\def\nothing{}%
\ifx\thearg\nothing
\vrule width0pt height\cellsize depth0pt\else
\hbox to 0pt{\usebox{\cell} \hss}\fi%
\vbox to \cellsize{
\vss
\hbox to \cellsize{\hss$#1$\hss}
\vss}}
\newcommand\tableau[1]{\vtop{\let\\\cr
\baselineskip -16000pt \lineskiplimit 16000pt \lineskip 0pt
\ialign{&\cellify{##}\cr#1\crcr}}}
%

\theoremstyle{plain}
\newtheorem{thm}{Theorem}
\newtheorem{lem}[thm]{Lemma}
\newtheorem{cor}[thm]{Corollary}
\newtheorem{prop}[thm]{Proposition}
\newtheorem{conj}[thm]{Conjecture}

\theoremstyle{definition}
\newtheorem{dfn}[thm]{Definition}
\newtheorem{exmp}[thm]{Example}
\newtheorem{remark}[thm]{Remark}
\newtheorem{fact}[thm]{Fact}



\newcommand{\thmref}[1]{Theorem~\ref{#1}}
\newcommand{\corref}[1]{Corollary~\ref{#1}}
\newcommand{\conjref}[1]{Conjecture~\ref{#1}}
\newcommand{\propref}[1]{Proposition~\ref{#1}}
\newcommand{\propsref}[1]{Propositions~\ref{#1}}
\newcommand{\lemref}[1]{Lemma~\ref{#1}}
\newcommand{\lemsref}[1]{Lemmas~\ref{#1}}
\newcommand{\factref}[1]{Fact~\ref{#1}}
\newcommand{\rmkref}[1]{Remark~\ref{#1}}
\newcommand{\remarkref}[1]{Remark~\ref{#1}}
\newcommand{\defref}[1]{Definition~\ref{#1}}
\newcommand{\egref}[1]{Example~\ref{#1}}
\newcommand{\figref}[1]{Figure~\ref{#1}}
\newcommand{\figsref}[2]{Figures~\ref{#1},\ref{#2}}
\newcommand{\secref}[1]{Section~\ref{#1}}

\newcommand{\cwd}{\operatorname{cwd}}
\newcommand{\rwd}{\operatorname{rwd}}
\newcommand{\ltx}{P_x}
\newcommand{\rtx}{Q_x}
\newcommand{\ltw}{P_w}
\newcommand{\rtw}{Q_w}
\newcommand{\rtq}{Q}
\newcommand{\rtqp}{Q'}
\newcommand{\rtqpp}{Q''}
\newcommand{\bigvert}{!{\vrule width 1pt}}

\newcommand{\miniprefig}
      {\scalebox{.25}{\includegraphics{hexpics/ministr.eps}}}

\newcommand{\minimslfig}
      {\scalebox{.08}{\includegraphics{mslpics/sec3412.eps}}}

\newcommand{\miniinvfig}
      {\scalebox{.08}{\includegraphics{mslpics/sec4231.eps}}}

\newcommand{\minihexafig}
      {\scalebox{.07}{\includegraphics{mslpics/minihexa.eps}}}
\newcommand{\minihexbfig}
      {\scalebox{.07}{\includegraphics{mslpics/minihexb.eps}}}
\newcommand{\minihexcfig}
      {\scalebox{.07}{\includegraphics{mslpics/minihexc.eps}}}
\newcommand{\minihexdfig}
      {\scalebox{.07}{\includegraphics{mslpics/minihexd.eps}}}
\newcommand{\minihexefig}
      {\scalebox{.07}{\includegraphics{mslpics/minihexe.eps}}}
\newcommand{\minihexffig}
      {\scalebox{.07}{\includegraphics{mslpics/minihexf.eps}}}
\newcommand{\minihexgfig}
      {\scalebox{.07}{\includegraphics{mslpics/minihexg.eps}}}

\newcommand{\mymslfig}[2]{\begin{figure}[htbp]\begin{center}
      {\scalebox{.4}{\includegraphics{#1}}}
      \caption{#2}\label{fig:#1}
    \end{center}\end{figure}}

\newcommand{\myintromslfig}[2]{\begin{figure}[htbp]\begin{center}
      {\scalebox{.4}{\includegraphics{#1.eps}}}
      \footnotesize{\textbf{Figure 0.1.} #2}
    \end{center}\end{figure}}

\newcommand{\mysmmslfig}[2]{\begin{figure}[ht]\begin{center}
      {\scalebox{.32}{\includegraphics{#1}}}
      \caption{#2}\label{fig:#1}
    \end{center}\end{figure}}

\newcommand{\mytinymslfig}[2]{\begin{figure}[ht]\begin{center}
      {\scalebox{.2}{\includegraphics{#1}}}
      \caption{#2}\label{fig:#1}
    \end{center}\end{figure}}

\newcommand{\myhexfig}[2]{\begin{figure}[ht]\begin{center}
    \input{hexpics/#1.pstex_t}\caption{#2}\label{fig:#1}
    \end{center}\end{figure}}

\newcommand{\mynewhexfig}[2]{\begin{figure}[ht]\begin{center}
    \input{hexpics/#1.pstex_t}\caption{#2}\label{fig:#1}
    \end{center}\end{figure}}

\newcommand{\mynewfig}[2]{\begin{figure}[ht]\begin{center}
      {\scalebox{.3}{\includegraphics{newpics/#1.eps}}}
      \caption{#2}\label{fig:#1}
    \end{center}\end{figure}}

\newcommand{\sfrac}[2]{\genfrac{\{}{\}}{0pt}{}{#1}{#2}}
\newcommand{\sumsb}[1]{\sum_{\substack{#1}}}  
\newcommand{\minsb}[1]{\substack{#1}}  
\newcommand{\pstack}[1]{P_{\minsb{#1}}}  
\newcommand{\musumr}[3]{\Theta_{(\cdot #3)}[#1,#2]}
\newcommand{\musuml}[3]{\Theta_{(#3 \cdot)}[#1,#2]}
\newcommand{\mscor}[3]{\delta_{(\cdot #3)}[#1,#2]}
\newcommand{\mslor}[3]{\omega_{(\cdot #3)}[#1,#2]}
\newcommand{\msset}[3]{\theta_{(\cdot #3)}[#1,#2]}
        
\newcommand{\mscol}[3]{\delta_{(#3 \cdot)}[#1,#2]}
\newcommand{\mslol}[3]{\omega_{(#3 \cdot)}[#1,#2]}

\newcommand{\fsk}{S_k}
\newcommand{\fsn}{S_n}
\newcommand{\fsp}[1]{S_{#1}}
\newcommand{\frg}{\mathfrak{g}}
\newcommand{\frh}{\mathfrak{h}}
\newcommand{\slgn}{\mathfrak{sl}_n}

\newcommand{\io}{{i_1}}
\newcommand{\iw}{{i_2}}
\newcommand{\inn}{{i_r}}

\newcommand{\bbA}{\mathbb{A}}
\newcommand{\bbR}{\mathbb{R}}
\newcommand{\bbC}{\mathbb{C}}
\newcommand{\bbZ}{\mathbb{Z}}
\newcommand{\bbN}{\mathbb{N}}
\newcommand{\bbQ}{\mathbb{Q}}

\newcommand{\wbi}{w^{\hat{i}}}
\newcommand{\xbi}{x^{\hat{i}}}
\newcommand{\zbi}{z^{\hat{i}}}
\newcommand{\wbio}{w^{\hat{1}}}
\newcommand{\wsbio}{ws^{\hat{1}}}
\newcommand{\xbio}{x^{\hat{1}}}
\newcommand{\zbio}{z^{\hat{1}}}

\newcommand{\tti}{t}
\newcommand{\wti}{{\widetilde{w}}}
\newcommand{\xti}{\protect{\widetilde{x}}}
\newcommand{\uti}{{\widetilde{u}}}
\newcommand{\vti}{{\widetilde{v}}}
\newcommand{\yti}{{\widetilde{y}}}
\newcommand{\zti}{{\widetilde{z}}}

\newcommand{\nts}{\negthickspace}

\newcommand{\trans}{\mathcal{T}}
\newcommand{\simref}{\mathcal{S}}
\newcommand{\calh}{\mathcal{H}}
\newcommand{\calw}{\mathcal{W}}


\newcommand{\msp}{\textsc{msp}}

\newcommand{\br}{\mathbf{a}}
\newcommand{\brp}{\mathbf{a'}}
\newcommand{\bsig}{{\boldsymbol{\sigma}}}
\newcommand{\bnu}{{\boldsymbol{\nu}}}
\newcommand{\bgam}{{\boldsymbol{\gamma}}}
\newcommand{\bdel}{{\boldsymbol{\delta}}}
\newcommand{\bal}{{\boldsymbol{\alpha}}}
\newcommand{\bbe}{{\boldsymbol{\beta}}}

\newcommand{\rds}{\operatorname{D_R}}
\newcommand{\lds}{\operatorname{D_L}}
\newcommand{\musum}{\operatorname{Musum}}
\newcommand{\rank}{\operatorname{lvl}}
\newcommand{\pt}{\operatorname{pt}}
\newcommand{\heap}{\operatorname{Heap}}
\newcommand{\lcz}{\operatorname{lcz}}
\newcommand{\rcz}{\operatorname{rcz}}
\newcommand{\mcz}{\operatorname{mcz}}
\newcommand{\ver}{\operatorname{ver}}
\newcommand{\codim}{\operatorname{codim}}
\newcommand{\bcone}{\operatorname{Cone_{\wedge}}}
\newcommand{\ucone}{\operatorname{Cone^{\vee}}}
\newcommand{\charp}{\operatorname{char}}

\newcommand{\IH}{\operatorname{IH}}
\newcommand{\ext}{\operatorname{Flush}}

\newcommand{\defeq}{:=}

\newcommand{\digx}{\operatorname{mat}(x)}
\newcommand{\digxt}{\operatorname{mat}(\widetilde{x})}
\newcommand{\digwt}{\operatorname{mat}(\widetilde{w})}
\newcommand{\augdigx}{\operatorname{mat}'(x)}
\newcommand{\digw}{\operatorname{mat}(w)}
\newcommand{\digp}[1]{\operatorname{mat}(#1)}
\newcommand{\dpw}[1]{d_{#1,w}}
\newcommand{\dpb}[2]{d_{#1,#2}}
\newcommand{\duv}{d_{u,v}}
\newcommand{\dvw}{d_{v,w}}
\newcommand{\dzw}{d_{z,w}}
\newcommand{\dxw}{d_{x,w}}
\newcommand{\dyw}{d_{y,w}}
\newcommand{\txw}{\Theta_{x,w}}
\newcommand{\dxtwt}{d_{\widetilde{x},\widetilde{w}}}

\newcommand{\rxjw}{\mathcal{R}(x_i,w)}
\newcommand{\rxjmw}{\mathcal{R}(x_{i-1},w)}
\newcommand{\rxow}{\mathcal{R}(x_0,w)}

\newcommand{\area}{\mathcal{A}}

\newcommand{\ph}{\phi}
\newcommand{\pht}{\phi_t}
\newcommand{\phs}{\phi_s}
\newcommand{\phtp}{\phi_{t'}}
\newcommand{\phtal}{\phi_{t_{\alpha_j}}}

\newcommand{\cP}{{\mathcal P}}
\newcommand{\xsing}{{X_w^{\text{sing}}}}
\newcommand{\maxsing}{\operatorname{\protect{Maxsing}}(X_\protect{w})}
\newcommand{\maxsingt}{\operatorname{Maxsing}(X_{\widetilde{w}})}
\newcommand{\schub}[1]{X_{#1}}

\newcommand{\sij}{s_{i_j}}
\newcommand{\sik}{s_{i_k}}
\newcommand{\sio}{s_{i_1}}
\newcommand{\sir}{s_{i_r}}


\newcommand{\puv}{P_{u,v}}
\newcommand{\pxw}{P_{x,w}}
\newcommand{\pww}{P_{w,w}}
\newcommand{\pyw}{P_{y,w}}
\newcommand{\pyv}{P_{y,v}}
\newcommand{\pzv}{P_{z,v}}
\newcommand{\pxws}{P_{x,ws}}
\newcommand{\pxsws}{P_{xs,ws}}
\newcommand{\pxtwt}{P_{\xti,\wti}}
\newcommand{\pzw}{P_{z,w}}
\newcommand{\pxiwi}{P_{x^{-1},w^{-1}}}
\newcommand{\psxw}{P_{sx,w}}
\newcommand{\pxsw}{P_{xs,w}}
\newcommand{\pxz}{P_{x,z}}
\newcommand{\pkl}[2]{P_{#1,#2}}


\newcommand{\delxw}{\Delta(x,w)}
\newcommand{\delyw}{\Delta(y,w)}
\newcommand{\delytw}{\Delta(yt,w)}

\newcommand{\cpw}{{\mathcal P}(\br)}
\newcommand{\cpww}{{\mathcal P}_w(\br)}
\newcommand{\cpws}{{\mathcal P}(\br s)}
\newcommand{\cpxw}{{\mathcal P}_x(\br)}
\newcommand{\cpyw}{{\mathcal P}_y(\br)}
\newcommand{\cpxsw}{{\mathcal P}_{xs}(\br)}
\newcommand{\cpxsws}{{\mathcal P}_{xs}(\br/s)}
\newcommand{\cpxws}{{\mathcal P}_x(\br/s)}
\newcommand{\cpzero}{{\mathcal P}_x^0(\br)}
\newcommand{\cpone}{{\mathcal P}_x^1(\br)}
\newcommand{\cpeps}{{\mathcal P}_x^\epsilon(\br)}

\newcommand{\cd}{{\mathcal D}}
\newcommand{\Dbr}{\Delta_\bsig}
\newcommand{\Dbrj}{\Delta_{\bsig[j]}}
\newcommand{\Dbrr}{\Delta_{\bsig[r]}}
\newcommand{\Dbrk}{\Delta_{\bsig[k]}}
\newcommand{\Dbrjm}{\Delta_{\bsig[j-1]}}

\newcommand{\apw}{P(\br)}
\newcommand{\apww}{P_w(\br)}
\newcommand{\apxw}{P_x(\br)}
\newcommand{\apxsw}{P_{xs}(\br)}
\newcommand{\apxws}{P_x(\br/s)}
\newcommand{\apxsws}{P_{xs}(\br/s)}

\newcommand{\apzw}{P_z(\br)}
\newcommand{\apzsw}{P_{zs}(\br)}
\newcommand{\apzws}{P_z(\br/s)}
\newcommand{\apzsws}{P_{zs}(\br/s)}

\newcommand{\apew}{P_e(\br)}
\newcommand{\cpew}{{\mathcal P}_e(\br)}


\newcommand{\fl}{\operatorname{{f}\,\!{l}}}
\newcommand{\ra}{\operatorname{unfl}}
\newcommand{\im}{\operatorname{Im}}
\newcommand{\slide}{\operatorname{\tau}}

\newcommand{\tpq}{t_{p,q}}
\newcommand{\tab}{t_{a,b}}
\newcommand{\talbt}{t_{\alpha,\beta}}
\newcommand{\tac}{t_{a,c}}
\newcommand{\tbc}{t_{b,c}}
\newcommand{\tcd}{t_{c,d}}
\newcommand{\tij}{t_{i,j}}
\newcommand{\tik}{t_{i,k}}
\newcommand{\tjk}{t_{j,k}}
\newcommand{\til}{t_{i,l}}
\newcommand{\talj}{t_{\alpha_j}}
\newcommand{\tgd}{t_{\gamma,\delta}}

\newcommand{\ttpq}{t_{p,q}}
\newcommand{\ttab}{t_{a,b}}
\newcommand{\ttalbt}{t_{\alpha,\beta}}
\newcommand{\ttac}{t_{a,c}}
\newcommand{\ttbc}{t_{b,c}}
\newcommand{\ttcd}{t_{c,d}}
\newcommand{\ttij}{t_{i,j}}
\newcommand{\ttik}{t_{i,k}}
\newcommand{\ttjk}{t_{j,k}}
\newcommand{\ttalj}{t_{\alpha_j}}
\newcommand{\ttgd}{t_{\gamma,\delta}}

\newcommand{\rp}[2]{\mathcal{R}({#1},{#2})}
\newcommand{\qp}[2]{\mathcal{E}_{{#1},{#2}}(x,w)}
\newcommand{\qpt}[2]{\mathcal{E}_{{#1},{#2}}(x,w)}

\newcommand{\rxw}{\mathcal{R}(\protect{x},\protect{w})}
\newcommand{\rytw}{\mathcal{R}(yt,w)}
\newcommand{\rysw}{\mathcal{R}(ys,w)}
\newcommand{\rysiw}{\mathcal{R}(ys_i,w)}
\newcommand{\qabxw}{\mathcal{E}_{a,b}(x,w)}
\newcommand{\qptxw}{\mathcal{E}_t(x,w)}
\newcommand{\qptabxw}{\mathcal{E}_{\tab}(x,w)}
\newcommand{\qbcxw}{\mathcal{E}_{b,c}(x,w)}
\newcommand{\qcdxw}{\mathcal{E}_{c,d}(x,w)}
\newcommand{\qbcyw}{\mathcal{E}_{b,c}(y,w)}

\newcommand{\qabxtwt}{\mathcal{E}_{a,b}(\widetilde{x},\widetilde{w})}
\newcommand{\qptxtwt}{\mathcal{E}_t(\widetilde{x},\widetilde{w})}
\newcommand{\qptabxtwt}{\mathcal{E}_{\tab}(\widetilde{x},\widetilde{w})}
\newcommand{\qbcxtwt}{\mathcal{E}_{b,c}(\widetilde{x},\widetilde{w})}
\newcommand{\qcdxtwt}{\mathcal{E}_{c,d}(\widetilde{x},\widetilde{w})}

\newcommand{\rxiwi}{\mathcal{R}(x^{-1},w^{-1})}
\newcommand{\rysiwi}{\mathcal{R}(sy^{-1},w^{-1})}
\newcommand{\ryw}{\mathcal{R}(y,w)}
\newcommand{\ryiwi}{\mathcal{R}(y^{-1},w^{-1})}
\newcommand{\rxtwt}{\mathcal{R}(\widetilde{x},\widetilde{w})}
\newcommand{\rxtw}{\mathcal{R}(xt,w)}

\newcommand{\ykm}{y_{k,m}}
\newcommand{\xkm}{x_{k,m}}
\newcommand{\xktm}{y_{k,m}}
\newcommand{\xrs}{x_{r,s}}
\newcommand{\xrts}{y_{r,s}}
\newcommand{\wkm}{w_{k,m}}
\newcommand{\vkm}{v_{k,m}}
\newcommand{\wktm}{v_{k,m}}
\newcommand{\wrs}{w_{r,s}}
\newcommand{\wrts}{v_{r,s}}
